\documentclass[a4paper, 12pt]{article}
\usepackage{amsmath, amssymb, mathrsfs,color}

\begin{document}
\newtheorem{thm}{Theorem}[section]
\newtheorem{lem}[thm]{Lemma}
\newtheorem{cor}[thm]{Corollary}
\newtheorem{conj}{Conjecture}
\newtheorem{rem}{Remark}[section]
\newtheorem{defn}{Definition}[section]
\allowdisplaybreaks[1]
\baselineskip 14pt

\title{Quasi-Invariance of Gaussian Measures \\ Transported by the Cubic NLS \\ with Third-Order Dispersion on $\mathbf{T}$}
\author{Arnaud DEBUSSCHE$^*$ and Yoshio TSUTSUMI$^{**}$}
\date{*Univ Rennes, CNRS, IRMAR - UMR 6625, F- 35000 Rennes, France. \\ **Department of Mathematics, Kyoto University, \\ Kyoto 606-8502, JAPAN}

\maketitle
\abstract{We consider the Nonlinear Schr\"odinger (NLS) equation and prove that the Gaussian measure  with covariance $(1-\partial_x^2)^{-\alpha}$ 
on $L^2(\mathbf T)$ is quasi-invariant for the associated flow for $\alpha>1/2$. This is sharp and improves a previous result obtained in \cite{OTT} where the values 
$\alpha>3/4$ were obtained. Also, our method is completely different and simpler, it is based on an explicit formula for the Radon-Nikodym derivative. We obtain an explicit formula for this latter in the same spirit as in  \cite{Cruz1} and \cite{Cruz2}. The arguments are general and can be used to other Hamiltonian equations.
}

\section{Introduction and Theorems}

In the present paper, we consider the transport of the Gaussian measures under the flow generated by the following Cauchy problem of the cubic NLS with third-order dispersion:
\begin{align}
   &\partial_t u = -i \bigl ( i \partial_x^3 + \beta \partial_x^2 \bigr ) u - i |u|^2u, \quad t \in \mathbf{R}, \ x \in \mathbf{T}, \label{3NLS} \\
   &u(0, x) = u_0(x), \quad x \in \mathbf{T},  \label{ic}
\end{align}
where $\beta$ is a real constant. It is known (see for instance \cite{MT1}) that this equation is globally well posed in $L^2(\mathbf{T})$.    

More precisely, we consider the Gaussian measure $\mu_\alpha$  
formally given by $\mu_\alpha(du) = C_\alpha e^{- \frac 1 2 \| u \|_{H^\alpha}^2} \, du$, where $C_\alpha$ is a normalization constant and study its evolution with respect to the flow of \eqref{3NLS}.
This measure $\mu_\alpha$ can be defined as  the distribution of the random variable $X = \sum_{n \in \mathbf{Z}} \frac {g_n(\omega)} {\langle n \rangle^\alpha} e^{inx}$ where $(g_n)_{n\in \mathbf{Z}}$ is a sequence of standard independent normal variable with values in $\mathbf{C}$. For $\alpha>1/2$, it can also be defined as the centered Gaussian measure on $L^2(\mathbf T)$ with covariance $(1-\partial_x^2)^{-\alpha}$ (for the definition of the Gaussian measure, see, e.g., \cite{K2} or \cite{DPZ} chapter 2).

For the last two decades, the probabilistic  approach to nonlinear evolution equations has attracted a lot of researchers in the field of PDEs as well as in the field of probability.
One of the interesting problems in this direction is to investigate how Gaussian measures are transported by such nonlinear Hamiltonian flows as Korteweg-de Vries, nonlinear Schr\" odinger equations and others.
Especially, it is natural to ask whether or not the Gaussian measure $\mu_\alpha$ is quasi-invariant, i.e., mutually absolutely continuous with respect to the transported Gaussian measure by the nonlinear Hamiltonian flow (see, e.g., \cite{Bour2}, \cite{FT}, \cite{GOTW}, \cite{OST}--\cite{PTV}, \cite{STX} and \cite{Tzv}).

From a viewpoint of probability theory, in \cite{R}, Ramer shows the Jacobi theorem for the transformation by general nonlinear mappings (see also Kuo \cite{K1}).
In \cite{Cruz1} and \cite{Cruz2}, Cruzeiro gives formulas for ordinary differential equations driven by vector fields in finite and infinite dimensions.
For nonlinear Hamiltonian dispersive equations, Tzvetkov \cite{Tzv} studied the quasi-invariance of Gaussian measures. He first uses Ramer's Theorem 
to obtain quasi-invariance for measure supported in high regularity Sobolev space and then develop a new method based on the evolution of the measure of
Borel sets to treat lower regularity. 
Note that he did not give an explicit formula for the Radon-Nikodym derivative of the transported Gaussian measure with respect to $\mu_\alpha$.

Nowadays, there are many papers about the quasi-invariance of Gaussian measures transported by various nonlinear dispersive equations (see \cite{FT} for the fractional NLS, \cite{GOTW} for the nonlinear wave equation, \cite{OST} and \cite{OT1} for the fourth order NLS and \cite{OT2} and \cite{PTV} for NLS).
In \cite{OTT}, Oh, the second author and Tzvetkov have showed the quasi-invariance of Gaussian measure $\mu_\alpha$ for $\alpha > 3/4$ under the flow generated by the Cauchy problem of the cubic NLS with third-order dispersion by using the Kuo-Ramer theorem (\cite{K1}, \cite{R}).

In this paper, we give a simpler proof of the result about the quasi-invariance of $\mu_\alpha$ under the flow of \eqref{3NLS}-\eqref{ic} given in \cite{OTT}  and bring down their assumption $\alpha > 3/4$ to $\alpha > 1/2$. This is optimal since for $\alpha\le 1/2$, we have $\mu_\alpha(H^s(\mathbf{T}))=0$ for any $s\ge 0$
and well posedness does not hold  in $H^s(\mathbf{T})$ for $s<0$ for \eqref{3NLS}. We treat only the case $\alpha \le 1$, since as explained in \cite{OTT} the case $\alpha >1$ is easy and does not need very refined arguments. 
We also give a more explicit formula of the Radon-Nikodym derivative than in \cite{Cruz2}. This formula is in fact general and generalizes to other hamiltonian systems. We believe that it has its own interest and can be used 
in other circumstances.  For instance (see Remark \ref{r3.3}), it can be used to simplify the proof of quasi-invariance for the fourth order NLS given in \cite{OST}.

Throughout this paper, we always assume that 
\begin{align}
    2 \beta/3 \not \in \mathbf{Z} \backslash \{ 0 \}. \tag{NR} \label{NR}
\end{align}
This is a kind of non-resonant condition (see \eqref{reson} and \eqref{nonreso} below). As already mentioned, for $s\ge 0$, equation \eqref{3NLS} is globally 
well posed in $H^s(\mathbf T)$. For $u_0\in H^s(\mathbf T)$, we denote by $u(t,u_0),\; t\in \mathbf R$,  the unique solution to \eqref{3NLS}.  For any measure 
$\nu$
supported by $L^2(\mathbf T)$, it makes sense to define the push-forward of this measure by the flow of \eqref{3NLS}. We say that 
$\nu$ is quasi-invariant if its push-forward is equivalent to $\nu$. 

  Our main theorems in this paper are the followoing.

\begin{thm} \label{qi}
Assume that \eqref{NR} is satisfied and $\alpha \in (1/2,1]$.
Then, for all $R>0$, the Gaussian measure $\chi_{\| u_0 \|_{\{L^2(\mathbf{T})} \leq R \}}\mu_\alpha(du_0)$ with $L^2$ norm cut-off weight is quasi-invariant under the flow generated by the third-order cubic NLS \eqref{3NLS}.
\end{thm}

The proof of this result is based on our second main result. To state this latter, we introduce
the following Cauchy problem of the third-order NLS truncated to the $(2N+1)$-dimensional space. 
\begin{align}
   &\partial_t u_N = -i \bigl ( i \partial_x^3 + \beta \partial_x^2 \bigr ) u_N - i P_N \bigl (|u_N|^2u_N \bigr ), \quad t \in \mathbf{R}, \ x \in \mathbf{T}, \label{3NLS0} \\
   &u_N(0, x) = u_{0,N}(x), \quad x \in \mathbf{T}.  \label{ic0}
\end{align}
where $P_N u = \frac 1 {\sqrt{2 \pi}} \sum_{|k| \leq N} (u, e^{ikx}) e^{-ikx}$. We will see that it has a unique solution $u_N(t, u_{0,N})$, defined for $t\in \mathbf R$. 
For $u_0\in H^s(\mathbf R)$, $s\ge 0$ we write $u_N(t, u_0)=u_N(t, P_Nu_0)$, $t\in \mathbf R$. 

\begin{thm} \label{rnd} Let $R>0$, $\alpha >1/2$ and $0 < s < \alpha -1/2$.
Let $u_0 \in H^s$.
Assume that the following functions 
\begin{align}
f_N(t,u_0)= \chi_{\{\| u_0 \|_{L^2(\mathbf{T})} \leq R \}}\exp\left( - \int_0^t \,\bigl ( i(|u_N|^2u_N)(-r,u_0), D^{2\alpha} u_N(-r,u_0) \bigr )  \ dr\right)\label{weight_N}
\end{align}
are uniformly bounded in $L^p(d\mu_\alpha)$ for some $p>1$. 
Then the following function
\begin{align}
f(t,u_0)=   \chi_{\{\| u_0 \|_{L^2(\mathbf{T})} \leq R \}}\exp\left( - \int_0^t \,\bigl ( i(|u|^2u)(-r,u_0), D^{2\alpha} u(-r,u_0) \bigr )  \ dr\right)\label{weight}
\end{align}
is in $L^p(d\mu_\alpha)$ and is the Radon-Nikodym derivative of  the Gaussian measure $\chi_{\{\| u_0 \|_{L^2(\mathbf{T})} \leq R \}}\mu_\alpha(du_0)$ with $L^2$ norm cut-off weight transported at time $t$  by \eqref{3NLS}.
\end{thm}

We have writen $(\cdot, \cdot) = \mathrm{Re} \, (\cdot, \cdot)_{L^2(\mathbf{T})}$ and denoted by $D$ the Fourier multiplier corresponding to $(1-\partial_{x}^2)^{1/2}$, see below. 
To prove Theorem \ref{qi}, it thus suffices to bound the $L^p(d\mu_\alpha)$ norm of  \eqref{weight_N}. 

\begin{rem} \label{comparison}
(i) Theorem \ref{qi} easily implies that for $\alpha\in(1/2,1]$, $\mu_\alpha$ is quasi-invariant. This is due to the preservation of the $L^2$ norm. As already mentioned, the result is known and easier to prove for $\alpha>1$ (see \cite{OTT}). 
 
(ii) In \cite{OTT}, it is showed that if (NR) is satisfied and $\alpha > 3/4$, the Gaussian measure $\mu_\alpha$ is quasi-invariant under the flow of \eqref{3NLS}-\eqref{ic}.
The proof in \cite{OTT} uses the Ramer theorem (see \cite{R}).
The assumptions in Ramer's theorem are formulated in terms of the linearized equation associated with the Cauchy problem \eqref{3NLS}-\eqref{ic} (see \cite[(1), (2) and (3) on page 166]{R}).
In the case of nonlinear dispersive equations, the linearized equation often does not work well, because it has less symmetry than the full system.
Indeed, in contrast to \cite{OTT}, we do not prove the smoothing effect of the nonlinear Duhamel term but that of the quadrilinear form associated with the nonlinear interaction of \eqref{3NLS} (see \eqref{nonlp} and \eqref{L2} below).
This makes our proof of Theorem \ref{qi} simpler and the lower bound of $\alpha$ smaller than in \cite{OTT}.

(iii) More generally, given $A$  a non negative self-adjoint operator on $L^2(\mathbf T)$ and  $V(u)$ be a smooth function of $u$ and $J = -i$, let us consider the 
equation:
\begin{align}
   &\frac {d u}{dt} = J \bigl ( A u + V'(u) \bigr ),  \quad t \in \mathbf{R},  \label{hamil}\\
   &u(0) = u_0,
\end{align}
and assume that it is well posed in $H^s(\mathbf T)$ for some $s$.  Consider the measure gaussian $\mu_\alpha$ of covariance $(1+A)^{-\alpha}$ for 
$\alpha$ such that $\mu_\alpha $ is supported by $H^s(\mathbf T)$.

Then, the following function is formally the Radon-Nikodym derivative  of the Gaussian measure $\chi_{\{\| u_0 \|_{L^2(\mathbf{T})} \leq R \}}\mu_\alpha(du_0)$ with $L^2$ norm cut-off weight transported by the flow is given by
\begin{align*}
   \chi_{\{\| u_0 \|_{L^2(\mathbf`{T})} \leq R \}}e^{\int_0^t \, \bigl ( JV'(u(-\sigma, u_0)), (1+A)^{\alpha} u(-\sigma, u_0) \bigr ) \, d\sigma}
\end{align*}
provided it defines an integrable function. The rigorous proof has to be done case by case for each equation but the proof we give below is quite general. Clearly, 
we can also consider more general domains than $\mathbf T$.

Assumptions of Theorem \ref{rnd} are clearly not optimal. For instance, it suffices to prove uniform integrability of the functions $f_N(t,\cdot)$. We think that
such refinement is not necessary.

(iv) Theorem \ref{rnd} is in the same spirit as the result in  Cruzeiro in \cite[Th\' eor\` eme 1.4.1 on page 208 and (4.3) on page 226]{Cruz2}.
There  Cruzeiro gives a formula of the Radon-Nikodym derivative for more general vector fields and considers the Wiener measure in dimension $1$, if we forget about boundary conditions, the Wiener measure corresponds to $\mu_1$. The result in \cite{Cruz2} holds under very restrictive assumptions 
on the vector field. Its divergence with respect to the Wiener measure (the operator $\delta$ in \cite{Cruz2}) has to be exponentially integrable. In the case considered here of Hamiltonian systems, this divergence simplifies considerably due to the preservation of the $L^2$ norm and to the fact that 
$\mu_\alpha$ is invariant for the linear equation. Also, our integrability condition is on trajectories and we can take advantage of smoothing for the Hamiltonian
flow. We thus obtain  a more explicit formula under weaker conditions.
Indeed, in our setting, conditions (i)-(iii) in Th\' eor\` eme 1.4.1 of \cite{Cruz2} are not satisfied.
In this respect, our Theorem \ref{rnd} is a refinement over Theorem 1.4.1 in \cite{Cruz2} for nonlinear wave and dispersive equations with Hamiltonian structure.
\end{rem}

Here we list the notation which will be used throughout this paper.
We denote the Fourier transform in the spatial variable of function $f(x)$ by $\hat f(k)$.
Let the Fourier transform in both time and spatial variables of $f(t,x)$ denote $\tilde f(\tau, k)$.
For $a \in \mathbf{R}$, we put $\langle a \rangle = (1+ a^2)^{1/2}$.
Let $D u = \mathcal{F}^{-1} [ \langle k \rangle \hat u(k) ]$.
For $s$, $b \in \mathbf{R}$, we define the Fourier restriction space $X^{s,b}$ and its norm as follows (see Bourgain\cite[(7.14) on page 211]{Bour1}).
\begin{align*}
      &\| u \|_{X^{s,b}} = \Bigl ( \sum_{k \in \mathbf{Z}} \int_{\mathbf{R}} \langle k \rangle^{2s} \langle \tau + k^3 - \beta k^2 \rangle^{2b} |\tilde u(\tau,k)|^2 \ d\tau \Bigr )^{1/2},  \\
      &X^{s,b} = \bigl \{ u \in \mathscr{S}'(\mathbf{R}^2)  \big | \ u(t,x+ 2\pi) = u(t,x), \ \| u \|_{X^{s,b}} < \infty \bigr \}.
\end{align*}
For $T > 0$, we also define the localized Fourier restriction space on $(-T, T)$ as follows.
\begin{align*}
      &\| u \|_{X^{s,b}_T} = \inf \bigl \{ \| v \|_{X^{s,b}} \big | \ v \in X^{s,b}, \ v(t) = u(t) \ \textrm{on} \ (-T, T)  \bigr \},  \\
      &X^{s,b}_T = \bigl \{ u \in \mathscr{D}'(\mathbf{R} \times \mathbf{T})  \big | \ \| u \|_{X^{s,b}_T} < \infty \bigr \}.
\end{align*}
The rest of the paper is organized as follows.
In Section 2, we prove several lemmas which are used for the proof of Theorem \ref{qi}.
In Section 3, we show Theorem \ref{qi} by using the results proved in Section 2 and Theorem \ref{rnd}.
In Section 4, we describe the proof of Theorem \ref{rnd}.


\section{Preliminaries}

\begin{lem}\label{gStr}
Let $1/4 > s > 0$and $p > 2$.
Then, we have
\begin{align}
  &\| f \|_{L^{4+ \frac {16s}{1-4s}}(\mathbf{R}; L^4(\mathbf{T}))} \leq C \| f \|_{X^{s, (\frac {1+2s}{3}) +}}, \label{gStr1} \\
  &\| f \|_{L^p(\mathbf{R} ; L^2(\mathbf{T}))} \leq C \| f \|_{X^{0, (1/2-1/p) +}}. \label{gStr2} 
\end{align}
\end{lem}

\noindent Proof.
We have the following two inequalities.
\begin{align}
   &\| f \|_{L^4(\mathbf{R} \times \mathbf{T})} \leq C \| f \|_{X^{0, 1/3 +}}, \label{Str} \\
   &\| f \|_{L^\infty(\mathbf{R} ; L^2(\mathbf{T}))} \leq C \| f \|_{X^{0, 1/2 +}}. \label{TimeSob} 
\end{align}
Inequalities \eqref{Str} and \eqref{TimeSob} are the Strichartz estimate and the Sobolev embedding in $t$, respectively (for the Strichartz estimate, see, e.g., \cite[Proposition 2.4 on page 1710]{MT1}).

We first have \eqref{gStr2} by \eqref{TimeSob} and interpolation.
We now show \eqref{gStr1}.
Let $\theta = 1 - 4s$ and $p$, $q$ defined by:
\[
   1/p = \theta/4 + (1-\theta)/\infty,\quad
   1/q = \theta/4 + (1-\theta)/2.
\]
Then $p=4+\frac{16s}{1-4s}$, $q=4-\frac{16s}{1+4s}$. Since 
$\theta/3 + (1-\theta)/2 = \frac {1+2s} {3}$,  the interpolation between \eqref{Str} and \eqref{TimeSob} yields
\begin{align}
   \| f \|_{L^{4+\frac{16s}{1-4s}}(\mathbf{R}; L^{4-\frac{16s}{1+4s}}(\mathbf{T}))} \leq C \| f \|_{X^{0, (\frac {1+2s}{3}) +}}.  \label{interp}
\end{align}
It remains to use the Sobolev embedding $H^{s, 4-\frac{16s}{1+4s}}(\mathbf{T}) \subset L^4(\mathbf{T})$ to obtain \eqref{gStr1}

$\square$

For $s \in \mathbf{R}$ and $p \geq 1$, let $B_p^s$ denote $B_{pp}^s(\mathbf{T})$,
where $B_{pq}^s(\mathbf{T})$ is the Besov space (for the definition of the Besov space, see, e.g., Triebel \cite{Tr}). Recall that $B_2^s=H^s(\mathbf T)$.

\begin{lem} \label{int}Let $\{ g_n \}$ be a sequence of independent equidistributed complex centered Gaussian random variables.
Let $\alpha > 0$, $ p \ge  2$ and $s \ge  0$ be such that $\alpha - 1 + 1/p > s$.
For $r \geq 1$ and $B > 0$, we set
\begin{align}
   F(\omega) :=  \chi_{\bigl \{ \bigl (\sum_{n \in \mathbf{Z}} |g_n(\omega)|^2/\langle n \rangle^{2\alpha} \bigr )^{1/2} < B \bigr \}} \mathrm{exp} \,  \Bigl ( \bigl \| \sum_{n \in \mathbf{Z}} \frac {g_n(\omega)} {\langle n \rangle^\alpha} e^{inx} \bigr \|_{B^s_p(\mathbf{T})}^r \Bigr ). \notag
\end{align}
Assume that  $B > 0$ for $r<\frac{4 \alpha p}{p-2+2ps}$ and $B$ is sufficiently small for $r=\frac{4 \alpha p}{p-2+2ps}$.
Then, $F( \omega ) \in L^1(d \omega)$.
\end{lem}

\begin{rem} \label{remint}
Lemma \ref{int} implies that if $\alpha \leq 1$ and $\alpha - \frac12 > s \ge 0$, then for any $C > 0$,
\[
   \chi_{\bigl \{ \bigl (\sum_{n \in \mathbf{Z}} |g_n(\omega)|^2/\langle n \rangle^{2\alpha} \bigr )^{1/2} < B \bigr \}} \mathrm{exp} \,  \Bigl ( C \bigl \| \sum_{n \in \mathbf{Z}} \frac {g_n(\omega)} {\langle n \rangle^\alpha} e^{inx} \bigr \|_{B^s_2(\mathbf{T})}^4 \Bigr ) \in L^1(d\omega),
\]
since in this case $4<\frac{2\alpha}s$.
\end{rem}

\begin{cor}\label{c1.2}
Let $Z_\alpha$ be a Gaussian random variable in $L^2(\mathbf{T})$ with covariance $D^{-2\alpha}$ for $\alpha > 1-1/p+s$.
Then, 
$$
\mathbb E \left( \exp\left( \|Z_\alpha \|_{B^s_p(\mathbf{T})}^r \right)\chi_{\|Z_\alpha \|_{L^2(\mathbf{T})}<B}\right) <\infty,
$$
for $p,\, r,\, s, \, B$ as above. Equivalently, if $\mu_\alpha$ is a Gaussian measure  with covariance  $D^{-2\alpha}$ then
$$
\int  \exp\left( \|x \|_{B^s_p(\mathbf{T})}^r \right)\chi_{\|x \|_{L^2(\mathbf{T})}<B} \; d\mu_\alpha(x) <\infty.
$$
\end{cor}

\textit{Proof of Corollary  \ref{c1.2}}.
For any $h\in L^2(\mathbf{T})$:
$$\mathbb E \left(\left( Z_\alpha, h\right)_{L^2(\mathbf{T})}^2\right)=\|D^{-\alpha} h\|_{L^2(\mathbf{T})}^2. $$ 
Let $Y=D^{\alpha}Z$, then $Y$ is Gaussian 
and since:
$$
\mathbb E \left(\left( Y, h\right)_{L^2(\mathbf{T})}^2\right)=\| h\|_{L^2(\mathbf{T})}^2 
$$
it has covariance $Id$ and it is a white noise. Therefore, it can be identified with $\{ g_n \}$ be a sequence of independent equidistributed complex Gaussian random variable.
We deduce:
\begin{align*}
   &\mathbb E \left( \exp\left( \|Z_\alpha \|_{B^s_p(\mathbf{T})} \right)\chi_{\|Z_\alpha \|_{L^2(\mathbf{T})}<B}\right) \\
   &= \mathbb E \left( \exp\left( \|D^{-\alpha}  Y\|_{B^s_p(\mathbf{T})} \right)\chi_{\|D^{-\alpha} Y\|_{L^2(\mathbf{T})}<B}\right).
\end{align*}
This is precisely the integral of $F$.
$\square$

\bigskip


\textit{Proof of Lemma \ref{int}}. The case $p=2$ and $s=0$ is clear, we assume that either $p>2$ or $s>0$.
For simplicity, we may assume $g_n \sim \mathcal{N}(0, 1)$.
We follow the argument by Bourgain \cite[the proof of Lemma 3.10]{Bour2}.
We estimate the probability of the following event.
\begin{align}
   \mathbb{P} \Bigl [ \omega \big | \ \bigl \| \sum_{n \in \mathbf{Z}} \frac {g_n(\omega)} {\langle n \rangle^\alpha}e^{inx} \bigr \|_{B^s_p(\mathbf{T})} > \lambda, \ \bigl ( \sum_{n \in \mathbf{Z}} \frac {|g_n(\omega)|^2}{\langle n \rangle^{2\alpha}} \bigr )^{1/2} < B \Bigr ]. \label{prob}
\end{align}

We first decompose the sum $\sum_{n \in \mathbf{Z}} a_n$ dyadically into
\[
   \sum_{|n| \leq M} a_n + \sum_{N > M} \sum_{|n| \sim N} a_n,
\]
where $M = \bigl ( \frac \lambda {2B} \bigr )^{\frac 1 { \frac 1 2 - \frac 1 p + s}}$.
We have by the Sobolev embedding
\begin{align}
   \bigl \| \sum_{|n| \leq M} \frac {a_n} {\langle n \rangle^{\alpha}} e^{inx} \bigr \|_{B^s_p(\mathbf{T})} &\leq M^{ \frac 1 2 - \frac 1 p + s} \bigl \| \sum_{|n| \leq M} \frac {a_n} {\langle n \rangle^{\alpha}} e^{inx} \bigr \|_{L^2(\mathbf{T})}  \notag \\
   &\leq M^{ \frac 1 2 - \frac 1 p + s} B = \frac \lambda 2. 
\end{align}
Let $\{ \sigma_N \}_{N \geq M}$ be a sequence such that
\[
   \sigma_N \geq 0, \quad 0 < \sum_{N > M} \sigma_N \leq 1.
\]
If $\| \sum_{|n| \sim N} \frac {g_n}{\langle n \rangle^\alpha} e^{i n x} \|_{B^s_p(\mathbf{T})} \leq \frac 1 2 \sigma_n \lambda$ for all $N > M$, then
\[
   \sum_{N > M} \bigl \| \sum_{|n| \sim N} \frac {g_n} {\langle n \rangle^\alpha} e^{inx} \bigr \|_{B^s_p(\mathbf{T})} \leq \frac \lambda 2
\]
and
\[
   \bigl \| \sum_{n \in \mathbf{Z}} \frac {g_n} {\langle n \rangle^\alpha} e^{i n x} \bigr \|_{B^s_p(\mathbf{T})} \leq \frac \lambda 2.
\]
Therefore, we deduce that
\begin{align*}
   &\mathbb{P} \Bigl [ \omega | \ \bigl \| \sum_{n \in \mathbf{Z}} \frac {g_n(\omega)}{\langle n \rangle^\alpha} e^{i n x} \bigr \|_{B^s_p(\mathbf{T})}>\lambda, \ \bigl ( \sum_{n \in \mathbf{Z}} \frac {|g_n(\omega)|^2} { \langle n \rangle^{2\alpha}} \bigr )^{1/2} < B \Bigr ] \\
   &\leq \mathbb{P} \Bigl ( \| \sum_{|n| \sim N} \frac {g_n} {\langle n \rangle^\alpha} e^{i n x} \bigr \|_{B^s_p(\mathbf{T})} > \frac 1 2 \sigma_N \lambda \ \text{for some} \ N > M \Bigr ) \\
   &\leq \sum_{N > M} \mathbb{P} \Bigl ( \bigl \| \sum_{|n| \sim N} \frac {g_n}{\langle n \rangle^\alpha} e^{i n x} \bigr \|_{B^s_p(\mathbf{T})} > \frac 1 2 \sigma_N \lambda \Bigr ).
\end{align*}
It suffices to estimate each of the above events.

For each dyadic block, we have by the Sobolev embedding
\begin{align}
   \bigl \| \sum_{|n| \sim N} a_n e^{ i n x} \bigr \|_{B^s_p(\mathbf{T})} \lesssim N^{1/2-1/p +s} \bigl \| \sum_{|n| \sim N} a_n e^{i n x} \bigr \|_{L^2(\mathbf{T})} \quad (2 \leq p \leq \infty).  \label{Sob}
\end{align}
We deduce that on each event on the right hand side above:
\begin{align*}
    \frac 1 2 \sigma_N \lambda &< \bigl \| \sum_{|n| \sim N} \frac{g_n} {\langle n \rangle^\alpha} e^{inx} \bigr \|_{B^s_p(\mathbf{T})} \leq N^{1/2-1/p+s} \bigl \| \sum_{|n| \sim N} \frac{g_n} {\langle n \rangle^\alpha } e^{inx} \bigr \|_{L^2(\mathbf{T})}  \\
    &\leq N^{-(\alpha-1/2+1/p-s)} \bigl \| \sum_{|n| \sim N} g_n e^{i n x} \bigr \|_{L^2(\mathbf{T})},
\end{align*}
which yields
\begin{align}
   \bigl ( \sum_{|n| \sim N} |g_n(\omega)|^2 \bigr )^{1/2} > \frac 1 2 \sigma_N N^{\alpha- 1/2+1/p-s} \lambda. \label{prob1}
\end{align}
Since the sequence $\{ g_n \}_{|n| \sim N}$ is the $N$-dimensional Gaussian random variable, we have only to evaluate $\displaystyle (2 \pi)^{-N/2} A_N \int_{r>\frac 1 2 \sigma_N N^{\alpha- 1/2+1/p-s} \lambda} r^{N-1} e^{-r^2/2} \, dr$ to estimate the probability of the event \eqref{prob1}, where $A_N$ is the area of the $(N-1)$-dimensional unit hyperball.
So, the probability of the event \eqref{prob1} is bounded by
\[
   C e^{C N \log (\sigma_N N^{\alpha-1/2+1/p-s} \lambda ) - c \sigma_N^2 N^{2\alpha-1+ 2/p-2s} \lambda^2}
\]
for some $C$, $c > 0$

Now we choose $\kappa>0$ such that $2\alpha-1+2/p-2s > 1+\kappa$ and take $\sigma_N = c N^{-\kappa} + (M/N)^{\alpha-1/2+1/p-s}$ for $N > M$ and obtain
$$
   \sum_{N > M}  C e^{C N \log (\sigma_N N^{\alpha-1/2+1/p-s} \lambda ) - c \sigma_N^2 N^{2\alpha-1+ 2/p-2s} \lambda^2} \leq e^{-c_0 M^{2\alpha-1+2/p} \lambda^2}. 
$$

Gathering the above:
\begin{align}
\label{tail}
   \mathbb{P} \Bigl [ \omega \big | \ &\bigl \| \sum_{n \in \mathbf{Z}} \frac {g_n(\omega)} {\langle n \rangle^\alpha}e^{inx} \bigr \|_{B^s_p(\mathbf{T})} > \lambda, \ \bigl ( \sum_{n \in \mathbf{Z}} \frac {|g_n(\omega)|^2}{\langle n \rangle^{2\alpha}} \bigr )^{1/2} < B \Bigr ]   \notag \\
   &\le e^{-c_0 M^{2\alpha-1+2/p-2s}\lambda^2} 
\end{align}

The power of the exponential on the right side of \eqref{tail}  must be greater than $\lambda^r$ and so it follows that
\begin{align}
   \lambda^r < C B^{- \frac {2\alpha-1 + 2/p-2s}{1/2-1/p+s}} \lambda^{\frac{2\alpha-1+2/p-2s}{1/2-1/p+s} + 2} \label{restr2}
\end{align}
for large $\lambda > 0$.
This inequality \eqref{restr2} is satisfied if $r<\frac{8 \alpha p}{p-2+2ps}$ or if $r=\frac{8 \alpha p}{p-2+ 2ps}$ and $B$ is sufficiently small.
$\square$

We finally have the following lemma concerning the global well-posedness of the Cauchy problem \eqref{3NLS}-\eqref{ic} and the convergence property of solutions of the truncated equations \eqref{3NLS0}, \eqref{ic0}.

\begin{lem} \label{gw}
Let $s \geq 0$ and let $u_0 \in H^s(\mathbf{T})$.
Then, there exist the unique global solutions $u, \, u_N \in C(\mathbf{R}; H^s(\mathbf{T}))$ of \eqref{3NLS}-\eqref{ic} and \eqref{3NLS0}-\eqref{ic0}, with $u_{0,N}=P_Nu_0$, respectively such that
\begin{align}
   \sup_{t \in [-T,T]} &\| u(t) \|_{H^s} + \| u \|_{X^{s,1/2}_T} \leq C \| u_0 \|_{H^s},    \quad T > 0, \label{bnd1} \\
   &\| u(t) \|_{L^2} = \| u_0 \|_{L^2},  \quad t \in \mathbf{R}, \label{bnd2} \\
    \sup_{t \in [-T,T]} &\| u_N(t) \|_{H^s} + \| u_N \|_{X^{s,1/2}_T} \leq C \| u_0 \|_{H^s},    \quad T > 0, \label{bnd10} \\
   &\| u_N(t) \|_{L^2} = \| P_N u_0 \|_{L^2},  \quad t \in \mathbf{R}, \label{bnd20}
\end{align}
where $C$ is a positive constant dependent only on $\| u_0 \|_{L^2}$ and $T$.
Furthermore, for any $T > 0$,
\begin{align}
    \sup_{t \in [-T,T]} \| u(t) - u_N(t) \|_{H^s} + \| u - u_N \|_{X_T^{s,1/2}} \ \longrightarrow \ 0 \quad (N \to \infty). \label{bnd3}
\end{align}
\end{lem}

The proof of Lemma \ref{gw} follows from the Strichartz estimate and  the contraction argument (see, e.g.,  \cite[Theorem 1.1 on page 1708]{MT1} for the global existence of solution and \cite[Lemma 2.27 on page 8]{Bour2} for the latter assertion of Lemma \ref{gw}).


\section{Proof of Theorem \ref{qi}}

We prove Theorem \ref{qi} using Theorem \ref{rnd}, whose proof is given below. We need to bound uniformly in $N$ the expression given in \eqref{weight_N} 
by an expression which it is integrable. This is the content of Lemma \ref{L1} below. The conclusion follows since we obtain that the transported truncated Gaussian measure has a strictly positive density with respect
to the truncated gaussian measure. 


We first note that
\begin{align}
&\int_0^T\bigl ( i(|u|^2u)(-r,\cdot), D^{2\alpha} u(-r,\cdot) \bigr ) dr\\
   &=\mathrm{Im} \int_0^T \bigl ( (|u|^2u) (-r,\cdot), D^{2\alpha} u(-r,\cdot) \bigr )  \ dr \label{nonlp} \\
   &= \mathrm{Im} \int_0^T \bigl ( (|u|^2 - \frac 1 \pi \| u \|_{L^2(\mathbf{T})}^2 ) u(-r,\cdot), D^{2\alpha} u(-r,\cdot) \bigr ) \ dr. \notag
\end{align}
Instead of the left hand side of \eqref{nonlp}, we estimate the right hand side of \eqref{nonlp}, because the deduction of the squared $L^2$ norm removes resonant frequencies.

For $k = k_1 - k_2 + k_3$, we define a phase function $\varPhi$ as follows
\begin{eqnarray}
   &\varPhi (k, k_1, k_2, k_3) = (\tau + k^3 - \beta k^2) - (\tau_1 + k_1^3 - \beta k_1^2)  \notag \\
   & - (\tau_2 + k_2^3 + \beta k_2^2) - (\tau_3 + k_3^3 - \beta k_3^2) \notag \\
   &= 3(k_1 - k_2)(-k_2 + k_3)(k_3 + k_1 - 2 \beta / 3), \label{reson}
\end{eqnarray}
where $\tau = \tau_1 + \tau_2 + \tau_3$.

\begin{lem} \label{L1}
Let $\alpha \in (1/2,1]$.
We put $s = \alpha-1/2-\varepsilon$ for sufficiently small $\varepsilon > 0$.
Then, 
\begin{align}
    \Bigl | \mathrm{Im}  \int_0^T \sum_{k \in \mathbf{Z}} &\sum_{\begin{subarray}{c} k = k_1-k_2+k_3 \\ (k_1-k_2)(-k_2+k_3) \neq 0 \end{subarray}} e^{-i r \varPhi(k, k_1, k_2, k_3)}  \notag \\
   &\times \hat v(-r,k_1) \bar {\hat v}(-r,k_2) \hat v(-r,k_3) \langle k \rangle^{2\alpha} \bar {\hat v}(-r,k) \ dr \Bigr | \notag \\
   &\leq C\bigl ( \|u_0\|_{L^2(\mathbf{T})}, T \bigr ) \| u_0 \|_{H^s}^3, \label{L2}
\end{align}
where $\hat v(t,k) = e^{it(k^3-\alpha k^2)} \hat u(t,k)$.
Furthermore, by $F(u)$, we denote the functional: $u \mapsto \mathbf{R}$ on the left side of \eqref{L2}, where $u$ is a solution of \eqref{3NLS}-\eqref{ic}.
Let $( u_N )_N$ be a sequence of solutions to \eqref{3NLS0}-\eqref{ic0} with $u_{0,N} = P_N u_0$.
Then, $F(u_N)$ satisfies the same bound \eqref{L2} and 
$F(u_N) \longrightarrow F(u)$ as $N \to \infty$.
\end{lem}

\begin{rem} \label{res1}
(i) In the summation of the left hand side of \eqref{L2}, the restriction $(k_1-k_2)(-k_2+k_3) \neq 0$ comes from the deduction of the squared $L^2$ norm on the right hand side of \eqref{nonlp}. The left hand side of \eqref{L2} is in fact equal to the right hand side of \eqref{nonlp}.

(ii) Lemma \ref{L1} implies the smoothing effect of the quadrilinear form associated with the nonlinear interaction of \eqref{3NLS}, since $s = \alpha - 1/2 - \varepsilon$.
The smoothing type estimate has been investigated by many authors for various nonlinear dispersive equations (see, e.g., \cite{TT} for modified KdV, \cite{BIT} and \cite{ET1} for KdV, \cite{ET2} and \cite {PTV} for NLS and \cite{MT2} for equation \eqref{3NLS}).
The smoothing type estimate is also applied to other problems, for example, the unconditional uniqueness of solution (see, e.g., \cite{KO} and \cite{MPV} for modified KdV and \cite{GKO} for NLS).
\end{rem}

{\it Proof of Lemma \ref{L1}}.
We show inequality \eqref{L2}.
We define
\begin{eqnarray*}
   &M = \max \bigl \{ \bigl |\tau + k^3- \beta k^2 \bigr |, \bigl |\tau_1 + k_1^3- \beta k_1^2 \bigr |, \\
   &\bigl |\tau_2 + k_2^3 + \beta k_2^2 \bigr |, \bigl |\tau_3 + k_3^3- \beta k_3^2 \bigr | \bigr \}. 
\end{eqnarray*}
Assume that
\begin{align}
   |k_3 + k_1| \gtrsim |k|.  \label{max}
\end{align}
Since we have
\[
     \frac 1 2 |k_1-k_2| + \frac 1 2 |-k_2 + k_3| + \frac 1 2 |k_3+k_1| \geq |k_1 - k_2 + k_3| = |k|,
\] 
then, by the identity \eqref{reson}, for $k \in \mathbf{Z}$ with $|k|$ large, we have 
\begin{equation}
   \exists c > 0 ; \ M \geq \frac 1 4 |\varPhi| \geq c |k| |k_1 - k_2| |-k_2 + k_3|. \label{nonreso}
\end{equation}

We follow the proof of Lemma 2.4 in \cite{MT2} (see also \cite[\S 2.2. Nonlinear estimate: part 1]{OTT}).
Let $T$ be an arbitrarily fixed positive constant.
The function $f(x)$ is defined as follows.
\begin{align}
   \hat f(k) &:= -2 \textrm{Im} \int_{-T}^0  \sum^{\eqref{max}}_{\begin{subarray}{c} k = k_1-k_2+k_3 \\ (k_1-k_2)(-k_2+k_3) \neq 0 \end{subarray}} e^{i r \varPhi(k, k_1, k_2, k_3)}  \label{diff} \\
   &\times \hat v(r,k_1) \bar {\hat v}(r,k_2) \hat v(r,k_3) \langle k \rangle^{2\alpha} \bar {\hat v}(r,k) \ dr,  \notag
\end{align}
where $\sum^{\eqref{max}}$ denotes the sum over the region such that \eqref{max} is satisfied. By symmetry, the left hand side of \eqref{L2} is  bounded 
by $\sum_{k\in \mathbf Z} |\hat f (k)|$ plus two other similar terms. 

Integration by parts yields
\begin{align}
   \hat f(k) = 2 \textrm{Im} \biggl [ &\sum^{\eqref{max}}_{\begin{subarray}{c} k = k_1-k_2+k_3 \\ (k_1-k_2)(-k_2+k_3) \neq 0 \end{subarray}} (i \varPhi)^{-1} \langle k \rangle^{2\alpha} \label{ibp} \\
   &\times \bigl ( \hat u(-T,k_1) \bar {\hat u}(-T,k_2) \hat u(-T,k_3) \bar {\hat u}(-T,k)   - \hat u_0(k_1) \bar {\hat u}_0(k_2) \hat u_0(k_3) \bar {\hat u}_0(k) \bigl ) \notag \\
   + \int_{-T}^0 &\sum^{\eqref{max}}_{\begin{subarray}{c} k = k_1-k_2+k_3 \\ (k_1-k_2)(-k_2+k_3) \neq 0 \end{subarray}}  (i \varPhi)^{-1} \langle k \rangle^{2\alpha} \notag \\
   \times  &\Bigl ( \sum_{k_1 = k_{11}-k_{12}+k_{13}} \hat u(r,k_{11}) \bar {\hat u}(r,k_{12}) \hat u(r,k_{13}) \Bigr ) \notag \\
   &\times \bar {\hat u}(r,k_2) \hat u(r,k_3) \bar {\hat u}(r,k) \ dr \biggr ] + \textrm{other similar terms}. \notag
\end{align}
For the solution $u$, we define a function $\mathbf u$ such that $\mathbf u(r) = u(r)$ on $[-T, 0]$ and $\mathbf u(r) = 0$ $(r < -T$ or $r > 0)$.
For simplicity, we also denote the function $\mathbf u$ by $u$.
Let $s = \alpha - 1/2-\varepsilon$ for sufficiently small $\varepsilon > 0$.
We set
\begin{align*}
   A =  \sum_{k \in \mathbf{Z}} &\langle k \rangle^{2 \alpha} \Bigl | \sum^{\eqref{max}}_{\begin{subarray}{c} k = k_1-k_2+k_3 \\ (k_1-k_2)(-k_2+k_3) \neq 0 \end{subarray}} (i \varPhi)^{-1} \\
   &\times \hat u(-T,k_1) \bar {\hat u}(-T,k_2) \hat u(-T,k_3) \bar {\hat u}(-T,k) \Bigr |, \\
   B = \sum_{k \in \mathbf{Z}} &\langle k \rangle^{2\alpha} \Bigl | \int_{-T}^0 \sum^{\eqref{max}}_{\begin{subarray}{c} k = k_1+k_2+k_3 \\ (k_1+k_2)(k_2+k_3) \neq 0 \end{subarray}}  (i \varPhi)^{-1} \\
   &\times  \sum_{k_1 = k_{11}+k_{12}+k_{13}} \hat u(r,k_{11}) \bar {\hat u}(r,k_{12}) \hat u(r,k_{13})   \\
   &\times \bar {\hat u}(r,k_2) \hat u(r,k_3) \bar {\hat u}(r,k) \ dr \Bigr |. 
\end{align*}
We only show the estimates of these two typical terms $A$ and $B$ for the proof of Lemma \ref{L1}, since the other terms can be similarly estimated.

We first estimate $A$.
We put $s = \alpha - 1/2-\varepsilon$ for sufficiently small $\varepsilon > 0$.
Let $\hat w(t,k) = \langle k \rangle^{s} |\hat u(t, k)|$.
We may assume that
$$
   |k| \gtrsim \max \{ |k_1|, |k_2|, |k_3| \}.
$$
Otherwise, more than one of $k_1$, $k_2$ and $k_3$ are much larger than $|k|$.
Suppose that $|k_1|$ and $|k_2|$ are much larger than $|k|$ and that $|k_3|$ has the same size as $|k|$ or the size less than $|k|$.
In this case, \eqref{reson} implies that $|\varPhi| \gtrsim |k_1| |k_2| |k_1-k_2| \gtrsim |k|^2 |k_1-k_2|$.
Suppose that the three of $|k_1|$, $|k_2|$ and $|k_3|$ are much larger than $|k|$.
In that case, we have
\[
   |k_1-k_2| \sim |k_3|, \quad |-k_2+k_3| \sim |k_1|, \quad |k_3+k_1-2\alpha/3| \sim |k_2|,
\]
since $k =k_1-k_2+k_3$.
Therefore, \eqref{reson} implies that $|\varPhi| \gtrsim |k_1| |k_2| |k_3| \gtrsim |k|^3$.
Accordingly, these cases are easier to treat.
By symmetry we may assume that $|k_1| \lesssim |k_2| \lesssim |k_3|$.
Now we have the following three cases.

(Case 1) \quad  $|k_3| \gg |k_2|, |k_1|$,

(Case 2)  \quad  $|k_3| \sim |k_2| \gg |k_1|$,

(Case 3) \quad  $|k_3| \sim |k_2| \sim |k_1|$.

\noindent We first consider Case 1.
We note that in Case 1, we have by \eqref{nonreso}
\begin{align}
      \exists c > 0 ; \ M \geq \frac 1 4 |\varPhi| \geq c |k|^2 |k_1 - k_2|. \label{nonreso1} 
\end{align}
We have by \eqref{nonreso1}, the change of variables $k_1' = k_1 - k_2$ and the Schwarz inequality
\begin{align*}
   A \leq C  \sum_k &\sum_{|k_1| \lesssim |k|} \langle k_1 - k_2 \rangle^{-1} \langle k \rangle^{-1+2\varepsilon} \\
   &\times \sum_{\begin{subarray}{c} k-k_1 = -k_2 + k_3 \\ |k_3| \lesssim |k| \end{subarray}} |\hat u(k_1)| |\hat u(k_2)| \, \hat w(k_3)  \hat w(k) \\
   \leq C \sum_{k_1'} &\langle k_1' \rangle^{-2+2\varepsilon} \sup_{k_1' \in \mathbf{Z}} \bigl ( \sum_{k_2} |\hat u(k_1'-k_2)| \,  |\hat u(k_2)| \bigr ) \\
   &\times \sup_{k_1' \in \mathbf{Z}} \bigl ( \sum_{k} \hat w(k-k_1')\hat w(k) \bigr ) \\
   \leq C &\bigl ( \sum_{k_1'} \langle k_1' \rangle^{-2+2\varepsilon} \bigr ) \bigl ( \sum_k |\hat u(k)|^2 \bigr ) \bigl ( \sum_k \hat w(k)^2 \bigr )  \\
   \leq C &\|u(-T)\|_{L^2}^2 \| u(-T) \|_{H^s}^2 \leq C \bigl (\|u_0\|_{L^2(\mathbf{T})}, T \bigr ) \| u_0 \|_{H^s}^2,
\end{align*}
by Lemma \ref{gw}. We have omitted the dependance on $T$ in the first lines of the computation.

Case 2 is treated in the same way as above, since we have the following inequality similar to \eqref{nonreso1}.
\[
      \exists c > 0 ; \ M \geq \frac 1 4 |\varPhi| \geq c |k|^2 |-k_2 + k_3|.
\]
Now we consider Case 3.
In this case, the worst subcase is that two of $|k_1-k_2|$, $|-k_2 + k_3|$ and $|k_3 + k_1|$ are small and the other one is large.
For example, we suppose that $|k_3+ k_1| \ \gg \ |k_1-k_2|, \ |-k_2+k_3|$.
Then, we have
\begin{align}
      \exists c > 0 ; \ M \geq \frac 1 4 |\varPhi| \geq c |k| |k_1 - k_2| |-k_2+k_3|.  \label{nonreso2} 
\end{align}
We choose $\varepsilon > 0$ such that $3\varepsilon \leq s$.
Let now $\hat z(t,k) = \langle k \rangle^{3\varepsilon} |\hat u(t, k)|$.
We have by \eqref{nonreso2}, the change of variables $k_1' = k_1 - k_2$ and the Schwarz inequality, omitting again the dependance on $T$:
\begin{align*}
   A \leq C  \sum_k &\sum_{|k_1| \lesssim |k|} \langle k_1 - k_2 \rangle^{-1-\varepsilon}  \\
   &\times \sum_{\begin{subarray}{c} k-k_1 = -k_2 + k_3 \\ |k_3| \lesssim |k| \end{subarray}} |\hat u(k_1)| \, \hat z(k_2) \hat w(k_3) \hat w(k) \\
   \leq C \sum_{k_1'} &\langle k_1' \rangle^{-1-\varepsilon} \sup_{k_1' \in \mathbf{Z}} \bigl ( \sum_{k_2} |\hat u(k_1'-k_2)| \, \hat z(k_2) \bigr ) \\
   &\times \sup_{k_1' \in \mathbf{Z}} \bigl ( \sum_{k} \hat w(k-k_1') \hat w(k) \bigr ) \\
   \leq C &\bigl ( \sum_{k_1'} \langle k_1' \rangle^{-1-\varepsilon} \bigr ) \bigl ( \sum_k |\hat u(k)|^2 \bigr )^{1/2} \bigl ( \sum_k |\hat z(k)|^2 \bigr )^{1/2} \bigl ( \sum_k |\hat w(k)|^2 \bigr ) \\
   \leq C &\|u(-T)\|_{L^2} \|u(-T)\|_{H^{s}}^3 \leq C \bigl (\|u_0\|_{L^2(\mathbf{T})}, T \bigr ) \| u_0 \|_{H^s(\mathbf{T})}^3.
\end{align*}

We next estimate $B$.
The integration interval $(-T, 0)$ in the integral of $B$ can be extended to $\mathbf{R}$, since the support of $u(r)$ is contained in $[-T, 0]$.
Let $\tilde v_1(\tau,k) = \langle k \rangle^{s} |\tilde u(\tau,k)|$ and $\tilde v_2(\tau,k) = \langle k \rangle^{-\varepsilon} | \tilde u(\tau,k)|$.
Furthermore, we may assume that
\[
     |k| \gtrsim \max \{ |k_1|, |k_2|, |k_3|, |k_{11}|, |k_{12}|, |k_{13}| \} 
\]
for almost the same reason as in the above proof for the estimate of $A$.
We begin with the proof for Case 1.
We have by the Plancherel theorem, inequality \eqref{nonreso1}, Lemma \ref{gStr} \eqref{gStr1} with $s = \varepsilon$, \eqref{gStr2} with $p = 1/(2\varepsilon)$ and the $L^2$ norm conservation
\begin{align*}
   B \leq C \int_{\mathbf{R}^5} &\sum_{k} \langle k \rangle^{-1+6\varepsilon } \bigl ( \sum_{\begin{subarray}{c} k = k_1 - k_2 + k_3 \\ k_1 = k_{11}-k_{12}+k_{13} \\ |k_{11}|, |k_{12}|, |k_{13}|, |k_2|, |k_3| \lesssim |k| \end{subarray}} \tilde v_2(\tau_1-\tau_2, k_{11}) \tilde v_2(\tau_2-\tau_3, k_{12}) \\
   &\times \tilde v_2(\tau_3-\tau_4, k_{13}) \bigr ) \tilde v_2 (\tau_4-\tau_5, k_2) \tilde v_1(\tau_5, k_3) \langle k_1-k_2 \rangle^{-1}  \\
   &\times \tilde v_1(\tau_1, k)  \bigr ) \ d\tau_1 d\tau_2 d\tau_3 d\tau_4 d\tau_5 \\
   \leq & C \int_{\mathbf{R}} \bigl ( (D^{-1} (v_2^4), D^{-1+6\varepsilon} (v_1^2) \bigr ) \ dr \\
   \leq &C \int_{\mathbf{R}} \| D^{-1} (v_2^4) \|_{L^2(\mathbf{T})} \| D^{-1+6\varepsilon} (v_1^2) \|_{L^2(\mathbf{T})} \ dr \\
   \leq &C \int_{\mathbf{R}} \| v_2^4 \|_{L^1(\mathbf{T})}   \| v_1^2 \|_{L^1(\mathbf{T})} \ dr  \\
   \leq &C \| D^\varepsilon v_2 \|_{L^{4/(1-4\varepsilon)}(\mathbf{R}; L^4(\mathbf{T}))}^4 \| v_1 \|_{L^{1/(2\varepsilon)} (\mathbf{R}; L^2(\mathbf{T}))}^2 \\
   \leq &C \| u \|_{X^{s,1/2}}^2 \| u \|_{X^{0, 1/2}}^4 \leq C \bigl (\|u_0\|_{L^2(\mathbf{T})}, T \bigr ) \|u_0\|_{H^s(\mathbf{T})}^2.
\end{align*}
Here, at the fourth and the last inequalities, we have used the Sobolev embedding and Lemma \ref{gw}, respectively.

We next consider Case 2.
In this case, we have
\[
         \exists c > 0 ; \ M \geq \frac 1 4 |\varPhi| \geq c |k|^2 |- k_2+k_3|,
\]
which is similar to \eqref{nonreso1}.
Therefore, the proof for Case 2 is the same as that for Case 1 and so we omit it.

We finally consider Case 3. 
In this case, the worst subcase is the following.
\begin{align}
      &\exists c > 0 ; \ M \geq \frac 1 4 |\varPhi| \geq c |k| |k_1 - k_2| |-k_2+k_3|,  \label{nearlyreso1} \\
      &|k_1-k_2|, \ |-k_2+k_3| \ll |k|.   \label{nearlyreso2}
\end{align}
Without loss of generality, we may assume that $|k_{11}| \sim |k|$.
Let $\varepsilon$ be a positive constant with $\alpha-1/2 > 9\varepsilon$.
We put $\tilde v_1(\tau,k) = \langle k \rangle^{s-\varepsilon} |\tilde u(\tau,k)|$, $\tilde v_2(\tau,k) = \langle k \rangle^{7\varepsilon} |\tilde u(\tau,k)|$ and $\tilde v_3(\tau,k) = \langle k \rangle^{-\varepsilon} | \tilde u(\tau,k)|$.
We first note that
\[
   \langle k_1-k_2 \rangle^{-1} \langle -k_2+k_3 \rangle^{-1} \leq \langle k_1-k_2 \rangle^{-2} + \langle -k_2+k_3 \rangle^{-2}.
\]
Therefore, in the same way as in Case 1, we obtain by the Plancherel theorem, inequality \eqref{nonreso1}, Lemma \ref{gStr} \eqref{gStr1} with $s = \varepsilon$, \eqref{gStr2} with $p = 1/(2\varepsilon)$ and the $L^2$ norm conservation
\begin{align*}
   B \leq C \int_{\mathbf{R}^5} &\sum_{k} \Bigl ( \sum_{\begin{subarray}{c} k = k_1 - k_2 + k_3 \\ k_1 = k_{11}-k_{12}+k_{13} \\ |k_{11}|, |k_{12}|, |k_{13}|, |k_2|, |k_3| \lesssim |k| \end{subarray}} \bigl ( \langle k_1-k_2 \rangle^{-2} + \langle -k_2+k_3 \rangle^{-2} \bigr )  \\
   &\times \tilde v_2(\tau_1-\tau_2, k_{11}) \tilde v_3(\tau_2-\tau_3, k_{12})  \tilde v_3(\tau_3-\tau_4, k_{13}) \Bigr )   \\
   &\times \tilde v_3 (\tau_4-\tau_5, k_2) \tilde v_1(\tau_5, k_3) \tilde v_1(\tau_1, k)  \bigr ) \ d\tau_1 d\tau_2 d\tau_3 d\tau_4 d\tau_5 \\
   \leq & C \bigl [ \int_{\mathbf{R}} \bigl ( D^{-2} (v_2 v_3^3),  (v_1^2) \bigr ) \ dr + \int_{\mathbf{R}} \bigl ( v_2 v_3^2 v_1, D^{-2} (v_1 v_3) \bigr ) \ dr \bigr ]  \\
   \leq &C \bigl [ \int_{\mathbf{R}} \| D^{-2} (v_2 v_3^3) \|_{L^\infty(\mathbf{T})} \| v_1^2 \|_{L^1(\mathbf{T})} \ dr \\
   &+ \int_{\mathbf{R}} \| v_2 v_3^2 v_1 \|_{L^1(\mathbf{T})} \| D^{-2} (v_1 v_3) \|_{L^\infty(\mathbf{T})} \ dr \bigr ] \\
   \leq &C \bigl [ \int_{\mathbf{R}} \| v_2 v_3^3 \|_{L^1(\mathbf{T})}   \| v_1 \|_{L^2(\mathbf{T})}^2 \ dr   \\
   &+ \int_{\mathbf{R}} \| v_2 \|_{L^4(\mathbf{T})} \| v_3 \|_{L^4(\mathbf{T})}^2 \| v_1 \|_{L^4(\mathbf{T})}   \| v_1 v_3 \|_{L^1(\mathbf{T})} \ dr \bigr ] \\
   \leq &C \bigl [ \| D^\varepsilon v_2 \|_{L^{4/(1-4\varepsilon)}(\mathbf{R}; L^4(\mathbf{T}))} \| D^\varepsilon v_3 \|_{L^{4/(1-4\varepsilon)}(\mathbf{R}; L^4(\mathbf{T}))}^3 \| v_1 \|_{L^{1/(2\varepsilon)} (\mathbf{R}; L^2(\mathbf{T}))}^2 \\
   + &\| D^\varepsilon v_2 \|_{L^{4/(1-4\varepsilon)}(\mathbf{R}; L^4(\mathbf{T}))} \| D^\varepsilon v_3 \|_{L^{4/(1-4\varepsilon)}(\mathbf{R}; L^4(\mathbf{T}))}^2 \| D^\varepsilon v_1 \|_{L^{4/(1-4\varepsilon)}(\mathbf{R}; L^4(\mathbf{T}))}   \\
   &\times \| v_1 \|_{L^{1/(2\varepsilon)} (\mathbf{R}; L^2(\mathbf{T}))} \| v_3 \|_{L^{1/(2\varepsilon)} (\mathbf{R}; L^2(\mathbf{T}))} \bigr ]   \\
   \leq &C \| u \|_{X^{s,1/2}}^3 \| u \|_{X^{0, 1/2}}^3  \leq C \bigl (\|u_0\|_{L^2(\mathbf{T})}, T \bigr ) \|u_0\|_{H^s(\mathbf{T})}^3.
\end{align*}
Here, at the fourth, the last but one and the last inequalities, we have used the Sobolev embedding, the fact that $s > 8\varepsilon$ and Lemma \ref{gw}, respectively.
Accordingly, the estimate of $B$ is completed for all Cases 1--3.
Thus, we have proved the estimates of $A$ and $B$, which completes the proof of  inequality \eqref{L2}.

The same arguments can be used for the solutions of  \eqref{3NLS0}-\eqref{ic0} and prove that $F(u_N)$ satisfies exactly the same bound as $F(u)$.
Finally, we have
\begin{align*}
   \bigl | F(u) - F(u_N) \bigr | \leq \sum_{j=1}^4 F_j(u, u_N).
\end{align*}
Here, we define $F_j$ as follows.
\begin{align*}
   F_j = \Bigl | \mathrm{Im}  \int_0^T \sum_{k \in \mathbf{Z}} &\sum_{\begin{subarray}{c} k = k_1-k_2+k_3 \\ (k_1-k_2)(-k_2+k_3) \neq 0 \end{subarray}} e^{-i r \varPhi(k, k_1, k_2, k_3)}  \notag \\
   &\times \hat w_{j1}(-r,k_1) \bar {\hat w}_{j2}(-r,k_2) \hat w_{j3}(-r,k_3) \langle k \rangle^{2\alpha} \bar {\hat w}_{j4}(-r,k) \ dr \Bigr |, \\
\end{align*}
where $\hat w_{jl}(t,k)$ denotes one of  the following three factors.
\begin{align*}
   &e^{it(k^3-\alpha k^2)} \hat u(t,k),   \\
   &e^{it(k^3-\alpha k^2)} \hat u_N(t,k), \\
   &e^{it(k^3-\alpha k^2)} (\hat u(t,k) - \hat u_N(t,k)),
\end{align*}
and the last factor appears only onece for each $F_j$.
In the same way as above, we can obtain the following estimate.
\begin{align*}
   |F_j(u, u_N)| \leq C &\bigl ( 1 + \|u\|_{C([-T, T]; H^s)} + \|u_N\|_{C([-T, T]; H^s)} \\
   + &\| u \|_{X_T^{s, 1/2}} + \| u_N \|_{X_T^{s, 1/2}} \bigr )^5  \\
   &\times \bigl ( \|u - u_N \|_{C([-T, T]; H^s)} + \| u - u_N \|_{X_T^{s, 1/2}} \bigr ).
\end{align*}
This inequality and \eqref{bnd3} in Lemma \ref{gw} imply $F(u_N) \longrightarrow F(u)$ $(N \to \infty)$.
$\square$

\begin{rem} \label{res2}
In the definition of $B$, the sum over $k_{11}-k_{12}+k_{13} = k_1$ includes frequencies with $k_{11}-k_{12}=0$ or $-k_{12}+k_{13}=0$, though the sum over $k_1-k_2+k_3 = k$ contains neither frequencies with $k_1-k_2=0$ nor frequences with $-k_2+k_3 = 0$.
This has no influence on the estimate of $B$, because we do not use the modulation identity \eqref{reson} with respect to $k_{11}$, $k_{12}$ and $k_{13}$ in the above-mentioned proof.
\end{rem}

\begin{rem}\label{r3.3}
In \cite{OST}, Oh , Sosoe and Tzvetkov show the quasi-invariance of Gaussian measure $\mu_\alpha$ transported by the forth order cubic NLS for $3/4 > \alpha > 1/2$.
Their result can be proved by the same argument as above, which gives a simpler proof.
In \cite{OST}, they use the infinite iteration of normal form reduction.
This requires that they infinitely many times repeat the integration by parts and verify the convergence of the resulting function series.
\end{rem}


\section{Proof of Theorem 1.2}

In this section, we describe the proof of Theorem 1.2. Since \eqref{3NLS0} is a finite dimensional ODE, we know that for all $t\in \mathbf R$ the mapping
$\Phi_{N,t}\,:\,u_{0,N}\mapsto u(t,u_{0,N})$ is a $C^\infty$ diffeomorphism. 
It follows that the transported  measure $\mu_{N,t,R}$ equals the push forward of
the truncated gaussian measure $\nu_{N,\alpha,R}=\chi_{\{\| u_{0,N} \|_{L^2(\mathbf{T})} \leq R \}}\mu_{N,\alpha}(du_{0,N})$ by this diffeomorphism which has a 
density with respect to this measure - here, $\mu_{N,\alpha}$ is  $2N+1$ dimensional marginal of $\mu_\alpha$. This density is given by the change of variable formula:
$$
f_N(t,u_{0,N})=\det\left|D\Phi_{N,t}(u_{0,N})\right|^{-1} G_{\alpha,N}(\Phi_{N,t}^{-1}(u_{0,N})) G_{\alpha,N}^{-1}(u_{0,N}) 
$$
where $G_{\alpha,N}$ is the density of $\mu_{N,\alpha}$ with respect to the $2N+1$ dimensional Lebesgue measure. Note that we use the conservation of the $L^2$ norm. The above formula shows that this density does not depend 
on $R$ and  is smooth with respect to $(t,u_{0,N})$. 

Classically, it satisfies an evolution equation and we find an explicit form for this density. More precisely, for any $t_0\in \mathbf R$, we may write 
for a smooth function $\varphi$:
\begin{align*}
&\frac{d}{dt}\int \varphi(u_{0,N})f_N(t,u_{0,N})d\nu_{N,\alpha,R}(u_{0,N})\big|_{t=t_0}\\
&=\frac{d}{dt}\int \varphi(u_N(t,u_{0,N})) d\nu_{N,\alpha,R}(u_{0,N})\big|_{t=t_0}\\
&=\frac{d}{dt}\int \varphi(u_N(t+t_0,u_{0,N})) d\nu_{N,\alpha,R}(u_{0,N})\big|_{t=0}\\
&=\frac{d}{dt}\int \varphi(u_N(t,u_{0,N})) f_N(t_0,u_{0,N})d\nu_{N,\alpha,R}(u_{0,N})\big|_{t=0}.
\end{align*}
On the other hand, by the chain rule:
\begin{align*}
&\frac{d}{dt}\int \varphi(u_N(t,u_{0,N})f_N(t_0,u_{0,N})d\nu_{N,\alpha,R}(u_{0,N})\big|_{t=0}\\
&= \int\left( \nabla_{u_{0,N}}\varphi(u_{0,N}),  -i \bigl ( i \partial_x^3 + \beta \partial_x^2 \bigr ) u_{0,N} - i P_N \bigl (|u_{0,N}|^2u_{0,N}\bigr )\right) f_N(t_0,u_{0,N})d\nu_{N,\alpha,R}(u_{0,N}).
\end{align*}
where $\nabla_{u_{0,N}}$ is gradient with respect to $u_{0,N}$. We have 
$$
div_{u_{0,N}}( -i \bigl ( i \partial_x^3 + \beta \partial_x^2 \bigr ) u_{0,N} - i P_N \bigl (|u_{0,N}|^2u_{0,N}\bigr ))=0. 
$$
This is a general property of Hamiltonian system, indeed we can write:
$$
 -i \bigl ( i \partial_x^3 + \beta \partial_x^2 \bigr ) u_{0,N} - i P_N \bigl (|u_{0,N}|^2u_{0,N}\bigr )=J\nabla_{u_{0,N}}H(u_{0,N})
$$
where $H$ is the energy of the third order NLS equation and $J$ is the antisymmetric operator 
corresponding to the multiplication by $i$.  Then it suffices 
to write
$$
div_{u_{0,N}}\left(J\nabla_{u_{0,N}}H_N(u_{0,N})\right)= \mbox{Tr} (J D^2H_N(u_{0,N})),
$$
where $D^2H_N$ is the Hessian of $H$. Since $D^2H_N$ is symmetric and $J$ antisymmetric, we have by the properties of the trace:
$$
\mbox{Tr} (J D^2H_N(u_{0,N})) = \mbox{Tr} (D^2H_N(u_{0,N}) J^T)= -\mbox{Tr}  (D^2H_N(u_{0,N}) J)= -\mbox{Tr}( J D^2H_N(u_{0,N}) )
$$
so that this quantity vanishes. 

Moreover
$$
(  D^{2\alpha} u_{0,N} ,  -i \bigl ( i \partial_x^3 + \beta \partial_x^2 \bigr ) u_{0,N} )=0,
$$
as easily seen by Parseval identity. 
Assume for the moment that $\varphi$ has compact support in the open ball 
of radius $R$, then by integration by parts  and the above two cancelations:
\begin{align*}
&\frac{d}{dt}\int \varphi(u_{0,N})f_N(t,u_{0,N})d\nu_{N,\alpha,R}(u_{0,N})\big|_{t=t_0}\\\\
&= \int \varphi(u_{0,N})\left(\nabla_{u_{0,N}} f_N(t_0,u_{0,N}) ,  i \bigl ( i \partial_x^3 + \beta \partial_x^2 \bigr ) u_{0,N} + i P_N \bigl (|u_{0,N}|^2u_{0,N}\bigr )\right)d\nu_{N,\alpha,R}(u_{0,N})\\
&- \int \varphi(u_{0,N}) \left(
D^{2\alpha} u_{0,N} ,   i P_N \bigl (|u_{0,N}|^2u_{0,N}\bigr )\right)f_N(t_0,u_{0,N})d\nu_{N,\alpha,R}(u_{0,N}).
\end{align*}

We deduce that $f_N(t,u_{0,N})$ satisfies the transport equation:
\begin{align*}
\frac{d}{dt}f_N(t,u_{0,N})&= \left( \nabla_{u_{0,N}}f_N(t,u_{0,N}) ,  i \bigl ( i \partial_x^3 + \beta \partial_x^2 \bigr ) u_{0,N} +i P_N \bigl (|u_{0,N}|^2u_{0,N}\bigr )\right)\\
&-\left(D^{2\alpha} u_{0,N} ,  i P_N \bigl (|u_{0,N}|^2u_{0,N}\bigr )\right)f_N(t,u_{0,N}).
\end{align*}
Note that since $\nabla_{u_{0,N}}f_N(t,u_{0,N}) $ and $D^{2\alpha} u_{0,N}$ belong to the range of the orthogonal projector $P_N$, we could forget it in the expression above.

Recalling that $f_N(0,\cdot)=1$, we have an explicit formula for the density:
$$
f_N(t,u_{0,N})=\exp\left(-\int_0^t  \left(D^{2\alpha} u_{N}(-r,u_{0,N}) ,  i  \bigl (|u_{N}|^2u_{N}\bigr )(-r,u_{0,N})\right)dr \right).
$$
We extend $f_N(t,\cdot)$ to $L^2(\mathbf T)$ by $f_N(t,u_0)=f_N(t,P_Nu_0)$. Under our assumption, the sequence $(f_N(t,\cdot))_N$ is bounded in 
$L^p(d\mu_\alpha)$. Up to a subsequence, it converges weakly in $L^p(d\mu_\alpha)$. By Lemmas \ref{L1},  for any $u_0$, $f_N(t,u_{0})$ converges
to
$$
f(t,u_0)=\exp\left(-\int_0^t  \left(D^{2\alpha} u(-r,u_{0}) ,  i  \bigl (|u_{}|^2u_{}\bigr )(-r,u_{0})\right)dr \right).
$$
We deduce that $f(t,\cdot)$ is the $L^p(d\mu_\alpha)$ weak limit of $f_N(t,\cdot)$.

Write now
\begin{align*}
&\int \varphi(P_Nu_{0}) f_N(t,P_Nu_0)\chi_{\{\| P_Nu_{0} \|_{L^2(\mathbf{T})} \leq R \}}\mu_{\alpha}(du_{0})\\
&=\int \varphi(u_{0,N})f_N(t,u_{0,N})d\nu_{N,\alpha,R}(u_{0,N})\\
&=\int \varphi(u_N(t,u_{0,N}))d\nu_{N,\alpha,R}(u_{0,N})\\
&=\int \varphi(u_N(t,P_Nu_{0}))\chi_{\{\| P_Nu_{0} \|_{L^2(\mathbf{T})} \leq R \}}\mu_{\alpha}(du_{0}).
\end{align*}
By Lemmas \ref{gw} and \ref{L1}, we may let $N\to \infty$ and obtain for $\varphi$ continuous and bounded on $L^2(\mathbf T)$:
$$
\int \varphi(u_{0})f(t,u_{0})\chi_{\{\| u_{0} \|_{L^2(\mathbf{T})} \leq R \}}\mu_{\alpha}(du_{0})=\int \varphi(u(t,u_{0}))
\chi_{\{\| u_{0} \|_{L^2(\mathbf{T})} \leq R \}}\mu_{\alpha}(du_{0}).
$$
The result follows.
\bigskip

{\bf Aknowledgment:} A. Debussche is partially supported by the French government thanks to the "Investissements d'Avenir" program ANR-11-LABX-0020-0, Labex Centre Henri Lebesgue.
Y.Tsutsumi is partially supported by JSPS KAKENHI Grant-in-Aid for Scientific Research (B) (17H02853).
Y.T. is also grateful to people in ENS de Rennes for their kind hospitality, since most part of this work was done while he stayed at ENS de Rennes.

\end{document}